\documentclass[a4paper,10pt,reqno]{amsart}

\usepackage{a4wide,latexsym,amssymb,amstext,amscd,xypic}
\usepackage[centertags]{amsmath}
\usepackage{amsfonts,amsthm}

\usepackage[latin1]{inputenc}

\theoremstyle{plain}
\begingroup
 \newtheorem{theorem}{Theorem}[section]
 \newtheorem{lemma}[theorem]{Lemma}
 \newtheorem{proposition}[theorem]{Proposition}
 \newtheorem{corollary}[theorem]{Corollary}
 
\endgroup

\theoremstyle{definition}
\begingroup
 \newtheorem{definition}[theorem]{Definition}
 \newtheorem{example}[theorem]{Example}
 
\endgroup
  
\theoremstyle{remark}
\begingroup
 \newtheorem{remark}[theorem]{Remark}
 
\endgroup

\newcommand{\qdr}[1]{\mbox{$\left [ #1 \right ]$}}		
\newcommand{\grf}[1]{\mbox{$\left \{ #1 \right \}$}}		

\def\aa{\alpha}

\def\ee{\varepsilon}
\def\ff{\phi}

\def\mm{\mu}

\def\pp{\pi}
\def\qq{\theta}

\def\ss{\sigma}

\def\Ff{\Phi}
\def\Gg{\Gamma}

\def\Qq{\Theta}
\def\Ss{\Sigma}

\def\Yy{\Psi}

\newcommand{\gra}[1]{\`{#1}}  
\newcommand{\ac}[1]{\'{#1}}   

\newcommand{\mc}[1]{\mathcal{#1}}			

\newcommand{\cat}[1]{\begin{bf}#1\end{bf}}

\newcommand{\kc}{{\mc{C}}}

\newcommand{\ke}{{\mc{E}}}
\newcommand{\kf}{{\mc{F}}}

\newcommand{\ko}{{\mc{O}}}
\newcommand{\kp}{{\mc{P}}}

\newcommand{\kt}{{\mc{T}}}

\def\CC{\mathbb C}     							
\def\PP{\mathbb P}    							
\def\RR{\mathbb R}   							
\def\ZZ{\mathbb Z}    							

\def\to{\rightarrow}							

\def\cont{{\, \subseteq \, }} 					

\def\glpiu{\widetilde{GL^+ (2, \RR)}}						

\newcommand{\os}[2]{\overset{#1}{#2}}				
\newcommand{\us}[2]{\underset{#1}{#2}}				
\newcommand{\mor}[1][]{\xrightarrow{#1}}

\def\sup{\mathop{\rm sup}\nolimits}

\def\deg{\mathop{\rm deg}\nolimits}

\def\rk{\mathop{\rm rk \,}}

\def\Hom{\mathop{\rm Hom \, }\nolimits}
\def\Aut{\mathop{\rm Aut \, }\nolimits}
\def\Auteq{\mathop{\rm Auteq \, }\nolimits}

\def\Ext{\mathop{\rm Ext \, }\nolimits}

\def\dim{\mathop{\rm dim \, }\nolimits}

\def\min{\mathop{\rm min \, }\nolimits}

\def\Tr{\mathop{\rm Tr \, }\nolimits}

\def\Stab{\mathop{\rm Stab \, }\nolimits}

\newcommand{\coh}{{\cat{Coh}}}

\title[Stability conditions on curves]{Stability conditions on curves}
\author{Emanuele Macr\gra{\i}}
\date{}

\address{Max-Planck-Institut f\"{u}r Mathematik, Vivatsgasse 7, 53111
Bonn, Germany}
\email{macri@mpim-bonn.mpg.de}



\begin{document}

\begin{abstract}
We study some examples of Bridgeland-Douglas stability conditions on triangulated categories. From one side we give a complete description of the stability manifolds for smooth projective curves of positive genus. From the other side we study stability conditions on triangulated categories generated by an exceptional collection. In the case of the projective line this leads to the connectedness and simply-connectedness of the stability manifold.
\end{abstract}

\maketitle

\section{Introduction}

The notion of stability condition on a triangulated category has been introduced by Bridgeland in \cite{Brid1}, following physical ideas of Douglas \cite{Doug}.

A stability condition on a triangulated category is given by abstracting the usual properties of $\mm$-stability for sheaves on projective varieties; one introduces the notion of slope, using a group homomorphism from the Grothendieck group $K(\cat{T})$ of the triangulated category $\cat{T}$ to $\CC$ and then requires that a stability condition has generalized Harder-Narasimhan filtrations and is compatible with the shift functor.

The main result of Bridgeland's paper is that these stability conditions can be described via a parameter space of stabilities. This space becomes a (possibly infinite-dimensional) manifold, called the stability manifold and denoted by $\Stab (\cat{T})$, if an extra condition (local finiteness) is assumed.

Bridgeland studied finitely dimensional slices of these spaces for the case of elliptic curves in \cite{Brid1}, but he left two cases: the curves of genus greater than one and $\PP^1$.
For the curve of genus greater than one there is a simple solution, applying a technical lemma of Gorodentsev, Kuleshov and Rudakov (\cite{GKR}, Lemma 7.2).

For the case of $\PP^1$, the situation is slightly more involved, since in $D(\PP^1)$ there are bounded $t$-structures whose heart is an abelian category of finite length (for example the one induced by the equivalence of $D(\PP^1)$ with the derived category of representations of the Kronecker quiver \cite{Bond}) and also there are degenerate stability conditions. Anyway, it is again possible to have an explicit description of $\Stab (\PP^1)$, using the classification of exceptional objects on $D(\PP^1)$.
For an application of these results to the study of stability manifolds of projective spaces and del Pezzo surfaces see \cite{emolo}.

\par\medskip

The plan of the paper is as follows. In Section \ref{sec:stabcond} we give a short summary on stability conditions on triangulated categories and we describe the stability manifolds for derived categories of smooth projective curves of positive genus. In Section \ref{sec:quiverexcobj} we study stability conditions associated to a complete exceptional collection on a triangulated category $\cat{T}$. After recalling some basic facts about quivers and exceptional objects, we show how to naturally associate to a complete exceptional collection with no negative homomorphisms between its objects a heart of a bounded $t$-structure and then a family of stability conditions having this one as heart. In this way we can define a collection of open connected subsets $\Ss \cont \Stab (\cat{T})$ of maximal dimension, parametrized by the orbits of the action of the braid group on exceptional collections (up to shifts). In Section \ref{sec:examples} we examine some topological properties of these subsets $\Ss$ for the derived category of the quiver $P_n$, with two vertices and $n$ arrow from the first vertex to the second one (for $n=2$ this is precisely the Kronecker quiver). In this case $\Ss (P_n)$ is unique, simply connected and coincides with all $\Stab (P_n)$.


\section{Stability conditions on triangulated categories}\label{sec:stabcond}

In this section we give a summary of Bridgeland's paper \cite{Brid1}. Let $\cat{T}$ be an essentially small triangulated category.

\begin{definition}\label{def:main}
A \emph{stability condition} $\sigma = (Z, \mc{P})$ on $\cat{T}$ consists of a group homomorphism $Z : K(\cat{T}) \to \CC$, called the \emph{central charge}, and strongly full, additive subcategories $\mc{P} (\ff) \cont \cat{T}$, $\ff \in \RR$. They keep the following compatibilities:
(1) for any nonzero object $E \in \mc{P} (\ff)$, $Z(E) \neq 0$ and $Z(E) / |Z(E)| = \exp (i \pp \ff)$;
(2) $\forall \ff \in \RR$, $\mc{P} (\ff + 1) = \mc{P} (\ff) [1]$;
(3) if $\ff_1 > \ff_2$ and $A_j \in \mc{P} (\ff_j)$, $j=1,2$, then $\Hom_{\cat{T}} (A_1 , A_2) = 0$;
(4) for any nonzero object $E \in \cat{T}$ there is a finite sequence of real numbers $\ff_1 > ... > \ff_n$ and a collection of triangles $E_{j-1} \to E_j \to A_j$
with $A_j \in \mc{P} (\ff_j)$, for all $j$, $E_0 = 0$ and $E_n = E$.
\end{definition}

The collection of exact triangles in Definition \ref{def:main} (4) is called the \emph{Harder-Narasimhan filtration} of $E$ (\emph{HN filtration} for short). Note that the HN filtration is unique up to isomorphisms. We write $\ff^+_\ss (E) := \ff_1$, $\ff^-_\ss (E) := \ff_n$, and $m_\ss (E) := \sum_j |Z(A_j)|$.
From the definition, each subcategory $\mc{P} (\ff)$ is extension-closed and abelian. Its nonzero objects are said to be \emph{semistable} of phase $\ff$ in $\ss$, and the minimal objects (classically called simple objects, i.e.\ objects without proper subobjects or quotients) are said to be \emph{stable}.

For any interval $I \cont \RR$, $\mc{P} (I)$ is defined to be the extension-closed subcategory of $\cat{T}$ generated by the subcategories $\mc{P} (\ff)$, for $\ff \in I$. Bridgeland proved that, for all $\ff \in \RR$, $\mc{P} ((\ff, \ff + 1])$ is the heart of a bounded $t$-structure on $\cat{T}$. The category $\mc{P} ((0, 1])$ is called the \emph{heart} of $\ss$.

\begin{remark}\label{rmk:tstruct}
Let $H := \grf{z \in \CC \, : \, z = |z| \exp (i \pi \ff), \, 0< \ff \leq 1}$. If $\cat{A} \cont \cat{T}$ is the heart of a bounded $t$-structure and moreover it is an abelian category of finite length (i.e.\ artinian and noetherian) with a finite number of minimal objects, then by \cite[Proposition 5.3]{Brid1}, a group homomorphism $Z: K(\cat{A}) \to \CC$ such that $Z (E) \in H$ for all minimal objects $E \in \cat{A}$ (such a group homomorphism is called a \emph{stability function}), extends to a unique stability condition on $\cat{T}$.
\end{remark}

\begin{lemma}\label{lem:comparison}
Let $(Z, \mc{P})$ be a stability condition on $\cat{T}$. Assume that $\cat{A}$ is a full abelian subcategory of $\mc{P} ((0,1])$ and the heart of a bounded $t$-structure on $\cat{T}$. Then $\cat{A} = \mc{P} ((0, 1])$.
\end{lemma}
\begin{proof}
By the definition of bounded $t$-structure (cfr., for example, \cite[Lemma 3.2]{Brid1}), if an object $E$ is in $\mc{P} ((0,1])$ and not in $\cat{A}$, then there is a nonzero morphism either $A[k] \to E$ or $E \to A [-k]$, with $k > 0$, $A \in \cat{A}$. But for all $j \in \ZZ$, $\cat{A} [j]$ is an abelian subcategory of $\mc{P} ((0,1]) [j]$, which leads to a contradiction.
\end{proof}

A stability condition is called \emph{locally-finite} if there exists some $\ee > 0$ such that each quasi-abelian subcategory  $\mc{P} ((\ff - \ee , \ff + \ee))$ is of finite length.
In this way $\mc{P} (\ff)$ has finite length so that every object in $\mc{P} (\ff)$ has a finite Jordan-H\"older filtration into stable factors of the same phase. The set of stability conditions which are locally-finite will be denoted by $\Stab (\cat{T})$.

By \cite[Proposition 8.1]{Brid1} there is a natural topology on $\Stab (\cat{T})$ defined by the generalized metric (i.e.\ it may be infinite)
\begin{equation}\label{eq:metric}
d(\ss_1, \ss_2) := \us{0 \neq E \in \cat{T}}{\sup} \grf{|\ff_{\ss_2}^+ (E) - \ff_{\ss_1}^+ (E)|, |\ff_{\ss_2}^- (E) - \ff_{\ss_1}^- (E)|, |\log \frac{m_{\ss_2} (E)}{m_{\ss_1} (E)}|} \in \qdr{0, \infty}.
\end{equation}
We call $\Stab (\cat{T})$ the \emph{stability manifold} associated to $\cat{T}$.

\begin{theorem}\label{thm:bridmain}\emph{\cite[Theorem 1.2]{Brid1}}
For each connected component $\Ss \cont \Stab (\cat{T})$ there is a linear subspace $V(\Ss) \cont (K(\cat{T}) \otimes \CC)^{\vee}$ with a well-defined linear topology such that the natural map $\mc{Z} : \Ss \to V(\Ss)$, which maps a stability condition $(Z, \mc{P})$ to its central charge $Z$, is a local homeomorphism. In particular, if  $K(\cat{T}) \otimes \CC$ is finite dimensional, then $\Ss$ is a finite dimensional complex manifold. 
\end{theorem}

For later use, we remind that $V (\Ss)$ is defined as the set of $W \in (K(\cat{T}) \otimes \CC)^{\vee}$ such that
$$\| W \|_{\ss} := \sup \grf{\frac{| W (E) |}{| Z (E) |} \, : \, E \ \mbox{is} \ \ss \text{-semistable}} < \infty,$$
where $\ss = (Z, \kp)$ is any stability condition in $\Ss$.

\begin{remark}\label{rmk:action} \emph{\cite[Lemma 8.2]{Brid1}}
The generalized metric space $\Stab (\cat{T})$ carries a right action of the group $\glpiu$, the universal cover of $GL^+ (2, \RR)$, and a left action by isometries of the group $\Auteq (\cat{T})$ of exact autoequivalences of $\cat{T}$. The second action is defined in the natural way. For the first action, let $(G,f) \in \glpiu$, with $G \in GL^+ (2,\RR)$ and $f: \RR \to \RR$ is an increasing map, periodic of period $1$ such that $G \exp (2i\pi \ff) / | G \exp (2i\pi \ff) | = \exp (2i \pi f(\ff))$, for all $\ff \in \RR$. Then $(G, f)$ maps $(Z, \kp) \in \Stab (\kt)$ to $(G^{-1} \circ Z, \kp \circ f)$.
\end{remark}


\subsection{Stability conditions on curves of positive genus}

Let $C$ be a smooth projective curve over $\CC$ of positive genus.
In this case the Grothendieck group of $D(C) := D^b (\coh (C))$ is not of finite rank as in the examples we will see in the rest of the paper. Hence the stability manifold $\Stab (D(C))$ is infinite dimensional. We can restrict to a finite dimensional (more geometric) slice of it.
According to \cite{Brid1}, we define $\Stab (C)$ as the finite dimensional submanifold of $\Stab (D(C))$ consisting of locally finite \emph{numerical stability conditions}, i.e.\ locally finite stability conditions whose central charge factorizes through the singular cohomology $H^* (C,\ZZ)$ of $C$, via the Chern character. An example of numerical stability condition on the bounded derived category of a curve is the stability condition induced by $\mu$-stability for sheaves \cite[Example 5.4]{Brid1}.

The fundamental ingredient for studying $\Stab (C)$ is this technical lemma \cite[Lemma 7.2]{GKR}.

\begin{lemma}\label{lem:GKR}
Let $C$ a smooth projective curve of genus $g(C) \geq 1$. Suppose $E \in \coh (C)$ is included in a triangle
$$Y \to E \to X \to Y \qdr{1}$$
with $\Hom^{\leq 0} (Y, X) = 0$. Then $X, Y \in \coh (C)$.  
\end{lemma}

\begin{theorem}
If $C$ has genus $g(C) \geq 1$, then the action of $\glpiu$ on $\Stab (C)$ is free and transitive, so that
$$\Stab (C) \cong \glpiu (\cong \CC \times {\mathbb H}),$$
where ${\mathbb H}$ denotes the complex upper half plane.
\end{theorem}
\begin{proof}
First note that the structure sheaves of points are stable. Indeed they are semistable because otherwise, by Lemma \ref{lem:GKR}, $\mc{O}_x$ is included in an exact sequence in $\coh (C)$
$$0 \to Y \to \mc{O}_x \to X \to 0,$$
with $\Hom^{\leq 0} (Y, X) = 0$, which is clearly impossible. Now, if $\mc{O}_x$ were not stable, then by the same argument all its stable factors should be isomorphic to a single object $K$, which implies $K \in \coh (C)$ and so $K \cong \mc{O}_x$.
In the same way, all line bundles are stable too.  

Then, let $\ss = (Z, \mc{P}) \in \Stab (C)$. Take a line bundle $A$ on $C$; by what we have seen above, $A$ and $\mc{O}_x$ are stable in $\ss$ with phases $\ff_A$ and $\psi$ respectively. The existence of maps $A \to \mc{O}_x$ and $\mc{O}_x \to A[1]$ gives inequalities $\psi -1 \leq \ff_A \leq \psi$, which implies that if $Z$ is an isomorphism (seen as a map from $H^* (C, \RR) \cong \RR^2$ to $\CC \cong \RR^2$) then it must be orientation preserving.
But $Z$ is an isomorphism: indeed if not, then there exist stable objects with the same phase having non-trivial morphisms, which is impossible. 
Hence, acting by an element of $\glpiu$, one can assume that $Z(E) = - \deg (E) + i \rk (E)$ and that for some $x \in C$, the skyscraper sheaf $\mc{O}_x$ has phase $1$. Then all line bundles on $C$ are stable in $\ss$ with phases in the interval $(0, 1)$ and all structure sheaves of points are stable of phase 1; but this implies that $\mc{P} ((0, 1]) = \coh (C)$ and so that the stability condition $\ss$ is precisely the one induced by $\mu$-stability on $C$.
\end{proof}


\section{Quivers and exceptional objects}\label{sec:quiverexcobj}

\subsection{Quivers and algebras}

In this subsection we give a quick review of some basic facts about finite dimensional algebras over $\CC$ and we start studying stability conditions on their derived categories.
For more details see \cite{Ben}.

A quiver is a directed graph, possibly with multiple arrows and loops. In this paper we deal only with finite quivers, that is those which have a finite number of vertices and arrows.
If $Q$ is a quiver, we define its path-algebra $\CC Q$ as follows. It is an algebra over $\CC$, which as a vector space has a basis consisting of the paths 
$$\bullet \to \ldots \to \bullet$$
in $Q$. Multiplication is given by composition of paths if the paths are composable in this way, and zero otherwise. Corresponding to each vertex $x$ there is a path of length zero giving rise to idempotent basis elements $e_x$.
Clearly $\CC Q$ is finitely generated as an algebra over $\CC$ if and only if $Q$ has only finitely many vertices and arrows, and finite dimensional as a vector space if and only if in addition it has no loops.

A representation of a quiver $Q$ associates to each vertex $x$ of $Q$ a vector space $V_x$, and to each arrow $x \to y$ a linear transformation $V_x \to V_y$ between the corresponding vector spaces. The dimension vector $\aa$ of such a representation is a vector of integers having length equal to the number of vertices of the quiver given by $\aa_x = \dim V_x$.
There is a natural one-to-one correspondence between representations of $Q$ and $\CC Q$-modules.
If $Q$ is finite without loops, then the simple modules correspond to vertices of the quivers (that is to representations which consist of a one dimensional vector space at a vertex and a zero dimensional vector space at any other vertex) and the indecomposable projective modules are of the form $P_x = \CC Q \cdot e_x$.

The importance of quivers in the theory of representations of finite dimensional algebras is illustrated by the following theorem (cfr. \cite[Proposition 4.1.7]{Ben}) due to Gabriel.

\begin{theorem} 
Every finite dimensional basic algebra over $\CC$ (i.e. every simple module is one dimensional over $\CC$) is the quotient of a path-algebra $\CC Q$ of a quiver $Q$ modulo an ideal contained in the ideal of paths of length at least two. In particular, if $\CC Q$ is finite dimensional, there is a bijection between the simple modules for the algebra and the simple $\CC Q$-modules.    
\end{theorem}

\begin{definition}
Let $I$ be a two-sided ideal contained in the ideal of paths of length at least two of the path-algebra of a quiver $Q$.
We call the pair $(Q, I)$ a quiver with relations (sometimes, with abuse of notation, we will forget the ideal $I$). The path-algebra of a quiver with relations is the algebra $\CC Q / I$. 
\end{definition}

As for quivers, we can define representations of a quiver with relations $Q$. Again there is a natural one-to-one correspondence between representations of $Q$ and modules over its path-algebra.

\begin{example}\label{ex:Pn}
The quiver $P_n$ ($n \geq 1$) contains two vertices and $n$ arrows from the first to the second vertex. For example,
$$P_2: \bullet \rightrightarrows \bullet.$$
\end{example}

\begin{example}\label{ex:SN}
The quiver with relations $T_N$ ($N \geq 1$) contains $N+1$ vertices $X_0, \ldots , X_N$ and $N (N+1)$ arrows $\ff_i^j : X_{i} \to X_{i+1}$ ($i=0, \ldots ,N-1$, $j = 0, \ldots ,N$). The relations are $\ff^j_{i+1} \ff^k_i = \ff^k_{i+1} \ff^j_i$.
Note that $T_1$ coincides with $P_2$.
\end{example}

In the next subsection we will see the connections between exceptional objects in derived categories and quivers with relations. We conclude now by examining stability conditions $(Z, \mc{P})$ on the derived category of a finite dimensional algebra $A$ for which $\mc{P} ((0, 1]) = \cat{mod}\text{-} A$, where $\cat{mod}\text{-} A$ denotes the category of finitely generated (right) $A$-modules.

\begin{lemma}\label{lem:fgalgebras}
Let $A$ be a finite-dimensional algebra over $\CC$ with simple modules $\grf{L_0 , \ldots, L_n}$ and let $(Z, \mc{P})$ be a stability condition on $D(A) := D^b (\cat{mod}\text{-} A)$. Assume that $L_0, \ldots L_n \in \mc{P} ((0, 1])$. Then $\cat{mod}\text{-} A = \mc{P} ((0, 1])$ and $L_j$ is stable, for all $j=0, \ldots, n$.
\end{lemma}
\begin{proof}
Since $\cat{mod}\text{-} A$ is the extension closed subcategory of $D(A)$ generated by $L_0, \ldots , L_n$, then $\cat{mod}\text{-} A$ is an abelian subcategory of $\mc{P} ((0, 1])$; so, the first part follows immediately from Lemma \ref{lem:comparison}.
Now the second statement is clear, since $L_j$ is a minimal object of $\text{mod-} A$.   
\end{proof}

\begin{remark}
Note that, as observed in \cite[Theorem 3]{Asp2},  in the situation of the previous lemma, an object of $D (A)$ is (semi)stable if and only if it is a shift of a $\qq$-(semi)stable $A$-module in the sense of King \cite{King}. In particular one can construct moduli spaces of semistable objects having fixed dimension vector.
\end{remark}


\subsection{Exceptional objects}

Some references for this subsection are \cite{Bond}, \cite{GR}, \cite{KO}.
Let $\cat{T}$ be a triangulated category linear over $\CC$ and of finite type, i.e. for any two objects $A, B \in \cat{T}$ the $\CC$-vector space $\oplus_{k \in \ZZ} \Hom^k (A, B)$ is finite-dimensional. 
Following \cite{Bond} we introduce the following notation for the graded complex of $\CC$ vector spaces with trivial differential
$$\Hom^\bullet (A, B) = \us{k \in \ZZ}{\bigoplus} \Hom^k (A, B) [-k],$$
where $A, B \in \cat{T}$, $\Hom^k (A, B) = \Hom (A, B [k])$. When $\cat{T}$ is the derived category of an abelian category, $\Hom^\bullet (A, B)$ is quasi-isomorphic to ${\mathbf R} \Hom (A, B)$.

\begin{definition}

(i) An object $E \in \cat{T}$ is called \emph{exceptional} if it satisfies
\begin{align*}
\Hom^i & (E, E) = 0, \qquad {\text{for}} \  i \neq 0,\\
\Hom^0 & (E, E) = \CC. 
\end{align*}

(ii) An ordered collection of exceptional objects $\grf{E_0, \ldots , E_n}$ is called \emph{exceptional} in $\cat{T}$ if it satisfies
$$\Hom^\bullet (E_i , E_j) = 0, \qquad {\text{for}} \  i > j.$$
We call an exceptional collection of two objects an \emph{exceptional pair}.
\end{definition}

\begin{definition}
Let $(E, F)$ an exceptional pair. We define objects $\mc{L}_E F$ and $\mc{R}_F E$ (which we call \emph{left mutation} and \emph{right mutation} respectively) with the aid of distinguished triangles
$$\mc{L}_E F \to \Hom^\bullet (E, F) \otimes E \to F,$$
$$E \to \Hom^\bullet (E, F)^* \otimes F \to \mc{R}_F E,$$
where $V [k] \otimes E$ (with $V$ vector space) denotes an object isomorphic to the direct sum of $\dim V$ copies of the object $E [k]$. Note that under duality of vector spaces the grading changes sign.
\end{definition}

A \emph{mutation} of an exceptional collection $\mc{E} = \grf{E_0, \ldots , E_n}$ is defined as a mutation of a pair of adjacent objects in this collection:
\begin{align*}
\mc{R}_i \mc{E} =& \grf{E_0, \ldots ,E_{i-1}, E_{i+1}, \mc{R}_{E_{i+1}} E_i, E_{i+2}, \ldots ,E_n},\\
\mc{L}_i \mc{E} =& \grf{E_0, \ldots ,E_{i-1}, \mc{L}_{E_i} E_{i+1}, E_i, E_{i+2}, \ldots ,E_n},
\end{align*}
$i = 0, \ldots , n-1$.
We can do mutations again in the mutated collection. We call any composition of mutations an \emph{iterated mutation}.

\begin{proposition}\label{pro:braidgroup}\emph{\cite{Bond}}
(i) A mutation of an exceptional collection is an exceptional collection.

(ii) If an exceptional collection generates $\cat{T}$ then the mutated collection also generates $\cat{T}$.

(iii) The following relations hold:
\begin{align*}
\mc{R}_i \mc{L}_i = \mc{L}_i \mc{R}_i = 1 && \mc{R}_{i} \mc{R}_{i+1} \mc{R}_{i} = \mc{R}_{i+1} \mc{R}_{i} \mc{R}_{i+1}  && \mc{L}_{i} \mc{L}_{i+1} \mc{L}_{i} = \mc{L}_{i+1} \mc{L}_{i} \mc{L}_{i+1}.
\end{align*}
\end{proposition}

The last relations, together with the obvious commutativity $\mc{R}_{i} \mc{R}_j = \mc{R}_j \mc{R}_i$ for $j-i \neq \pm 1$, could be rephrased by saying that there is an action of the braid group $A_{n+1}$ of $n+1$ strings on the set of exceptional collections.

\begin{definition}
Let $\mc{E} = \grf{E_0, \ldots , E_n}$ be an exceptional collection. We call $\mc{E}$
\begin{itemize}
\item strong, if $\Hom^k (E_i, E_j) = 0$ for all $i$ and $j$, with $k \neq 0$;
\item Ext, if $\Hom^{\leq 0} (E_i, E_j) = 0$ for all $i \neq j$;
\item regular, if $\Hom^k (E_i, E_j) \neq 0$ for at most one $k \geq 0$, for all $i$ and $j$;
\item orthogonal, if $\Hom^k (E_i, E_j) = 0$ for all $i \neq j$ and $k$;
\item complete, if $\mc{E}$ generates $\cat{T}$ by shifts and extensions.
\end{itemize}
\end{definition}

The relation between strong exceptional collections and finite dimensional algebras is contained in the following result due to Bondal.

\begin{theorem}\label{thm:Bondal} 
Let $\cat{T}$ be the bounded derived category of an abelian category. Assume that $\cat{T}$ is generated by a strong exceptional collection $\grf{E_0, \ldots , E_n}$. 
Then, if we set $E = \oplus E_i$ and $A = \Hom (E, E)$, $\cat{T}$ is equivalent to the bounded derived category of finite dimensional modules over the algebra $A$.   
\end{theorem}
\begin{proof}
Define the functor $\Ff : \cat{T} \to D(A)$ as the derived functor $\Ff (Y) = {\mathbf R} \Hom (E, Y)$, where $Y \in \cat{T}$, with the natural action of $A$ on that complex. For the proof that $\Ff$ is actually an equivalence see \cite[Theorem 6.2]{Bond}.  
\end{proof}

\begin{example}\label{ex:PN}\cite{Beil} \cite{GR}
Let $\cat{T} = D(\PP^N) := D^b (\coh (\PP^N))$. Then a complete strong exceptional collection is given by $\grf{\mc{O}, \ldots , \mc{O}(N)}$. The corresponding  algebra is given by the quiver with relations $T_N$ of Example \ref{ex:SN}.
\end{example} 

\begin{example}\label{ex:KO}\cite{KO} 
Let $\cat{T} = D(S) := D^b (\coh (S))$, where $S$ is a Del Pezzo surface, i.e. a smooth projective surface whose anticanonical class is ample. Then there exists a complete strong exceptional collection and all exceptional collections are obtained from this collection by iterated mutations.
\end{example}


\subsection{Stability conditions on triangulated categories generated by an exceptional collection}

In this subsection we study stability conditions on triangulated categories generated by Ext-exceptional collections.
Given a subcategory $S$ of $\cat{T}$, we denote by $\langle S \rangle$ the extension-closed subcategory of $\cat{T}$ generated by $S$, and by $\Tr (S)$ the minimal full triangulated subcategory containing $S$ and closed by isomorphisms.

\begin{lemma}\label{lem:tstruct}
Let $\grf{E_0, \ldots, E_n}$ be a complete Ext-exceptional collection on $\cat{T}$. Then $\langle E_0, \ldots, E_n \rangle$ is the heart of a bounded $t$-structure on $\cat{T}$.
\end{lemma}
\begin{proof}
We proceed by induction on $n$. If $n=0$ there is nothing to prove. Indeed $\cat{T} \cong D^b (\CC \text{-} \cat{vect})$, 
Assume $n >0$. Then consider the full triangulated subcategory $\Tr (E_1, \ldots, E_n)$ of $\cat{T}$.
This is an admissible subcategory \cite[Theorem 3.2]{Bond} and its right orthogonal is $\Tr (E_0)$. Moreover
$$\Tr (E_1, \ldots, E_n) \mor[i_*] \cat{T} \mor[j^*] \Tr (E_0)$$
is an exact triple of triangulated categories \cite[Proposition 1.6]{BK}.

By \cite[\S 1.4]{BBD} any pair of $t$-structures on $\Tr (E_1, \ldots, E_n)$ and $\Tr (E_0)$ determines a unique compatible $t$-structure on $\cat{T}$ given by
\begin{align*}
\cat{T}^{\leq 0} =& \grf{F \in \cat{T} \, : \, j^* F \in \Tr (E_0)^{\leq 0}, \ i^* F \in \Tr (E_1,\ldots ,E_n)^{\leq 0}}\\
\cat{T}^{\geq 0} =& \grf{F \in \cat{T} \, : \, j^* F \in \Tr (E_0)^{\geq 0}, \ i^! F \in \Tr (E_1,\ldots,E_n)^{\geq 0}},
\end{align*}
where $i^*$ and $i^!$ are respectively the left and right adjoint to $i_*$. More explicitly, if $F$ decomposes as
$$A \to F \to B$$
with $A \in \Tr (E_1, \ldots, E_n)$, $B \in \Tr (E_0)$, then $i^! F = A$; if $F$ decompose as
$$B' \to F \to A'$$
with $A' \in \Tr (E_1, \ldots, E_n)$, $B' \in \Tr (\mc{R}_{n-1} \ldots \mc{R}_0 E_0)$, then $i^* F = A'$.

By induction we can choose $t$-structures on $\Tr (E_0)$ and $\Tr (E_1,\ldots, E_n)$ having hearts $\langle E_0 \rangle$ and $\langle E_1 ,\ldots, E_n \rangle$ respectively.

We want to prove that $\langle E_0 , \ldots, E_n \rangle = \cat{T}^{\leq 0} \cap \cat{T}^{\geq 0} =: \cat{A}$.
Clearly $E_1, \ldots , E_n$ belong to $\cat{A}$. Moreover, by mutating $E_0$, since the exceptional collection is Ext, we have $i^* E_0 \in \Tr (E_1,\ldots ,E_n)^{\leq 0}$.
Hence $\langle E_0 ,\ldots, E_n \rangle$ is a full subcategory of $\cat{A}$.
If $F \in \cat{A}$, $F \notin \langle E_0 ,\ldots, E_n \rangle$, we can filter $F$ as
\begin{equation}\label{eq:tstruct}
A \to F \to B
\end{equation}
as before. By construction $B \in \langle E_0 \rangle$.
We want to prove that $A \in \langle E_1 ,\ldots, E_n \rangle$.
Assume the contrary. From the triangle \eqref{eq:tstruct} we get
$$0 \to H^0 (A) \to H^0 (F) \to H^0 (B) \to H^1 (A) \to 0.$$
Since $A \in \Tr (E_1 ,\ldots, E_n)$, by induction we can filter it as
$$C \to A \to D,$$
with $C \in \langle E_1 ,\ldots, E_n \rangle$ and $D \in \langle E_1 ,\ldots, E_n \rangle [-1]$.
But then $H^0 (A) \cong C$ and $H^1 (A) \cong D [1]$.
This means that we have a map $B \to D[1]$, a contradiction.

The fact that the glued $t$-structure is bounded is now clear.
\end{proof}

\begin{corollary}\label{cor:2objs}
Let $\mc{E} = \grf{E_0, \ldots, E_n}$ be a complete exceptional collection on $\cat{T}$ such that, for some $i < j$, $\Hom^{\leq 0} (E_i, E_j) =0$. Then $\langle E_i, E_j \rangle$ is a full abelian subcategory of $\cat{T}$.  
\end{corollary}

\begin{lemma}\label{lem:excobj}
Let $\grf{E_0, \ldots, E_n}$ be a complete Ext-exceptional collection on $\cat{T}$ and let $(Z, \mc{P})$ be a stability condition on $\cat{T}$. Assume $E_0, \ldots , E_n \in \mc{P} ((0, 1])$. Then $\cat{Q} := \langle E_0, \ldots , E_n \rangle =  \mc{P} ((0, 1])$ and $E_j$ is stable, for all $j=0, \ldots , n$.      
\end{lemma}
\begin{proof}
First of all notice that $E_0, \ldots , E_n$ are the minimal objects of $\cat{Q}$.
Indeed, let $0 \neq A \hookrightarrow E_i$ be a subobject in $\cat{Q}$. Then by definition there exists an inclusion $E_j \hookrightarrow A$, for some $j$. But, since $\grf{E_0, \ldots, E_n}$ is an Ext-exceptional collection, then $j$ should be equal to $i$. Hence the composite map $E_j \hookrightarrow A \hookrightarrow E_i$ is an isomorphism and $A \cong E_i$.
The proof is then the same as Lemma \ref{lem:fgalgebras}. 
\end{proof}

If two exceptional objects have at most one nontrivial $\Hom^k$, one can say more. 

\begin{proposition}\label{pro:mutstab}
Let $\grf{E_0, \ldots, E_n}$ be a complete Ext-exceptional collection on $\cat{T}$ and let $\ss = (Z, \mc{P})$ be a stability condition on $\cat{T}$ such that $E_0, \ldots, E_n \in \mc{P} ((0, 1])$. Fix $i < j$. Then $\ss$ induces a stability condition $\ss_{i j}$ on $\Tr (E_i, E_j)$ in such a way that every (semi)stable object in $\ss_{i j}$ with phase $\ff$ corresponds to a (semi)stable object in $\ss$ with the same phase. Moreover, assume that $\Hom^1 (E_i , E_j) \neq 0$ and $\Hom^k (E_i, E_j) =0$, for all $k \neq 1$. Then, if
$$\ff (E_j) \leq \ff (E_i),$$
$\mc{R}_{E_j} E_i$ and $\mc{L}_{E_i} E_j$ are semistable and, if $\ff (E_j) < \ff (E_i)$, they are stable.
\end{proposition}
\begin{proof}
Consider the triangulated category $\Tr (E_i , E_j)$ generated by $E_i$ and $E_j$.
Then $Z$ defines a stability function on the abelian category $\langle E_i, E_j \rangle$, and so a stability condition on $\Tr (E_i, E_j)$, by Remark \ref{rmk:tstruct}.

Let $S$ be a semistable object in $\Tr (E_i, E_j)$. We can assume $S \in \langle E_i, E_j \rangle$ and $\ff (E_j) \leq \ff (E_i)$.
Let
$$0 \to A \to S \to B \to 0$$
be a destabilizing filtration in $\cat{Q} := \langle  E_0, \ldots, E_n \rangle$. We have to prove that $A, B \in \langle E_i, E_j \rangle$.
If $A \in \langle E_{j_1},\ldots , E_{j_k} \rangle$, $j_1 < \ldots < j_k$ and $B \in \langle E_{i_1}, \ldots ,E_{i_s} \rangle$, $i_1 < \ldots < i_s$, then by hypothesis $j_k \leq j$, $i_1 \geq i$.
If $A'$ is the subobject of $A$ belonging to $\langle E_j \rangle$, then we have the diagram
\begin{equation*}
\xymatrix{& & & 0 \ar[d] \\
	        & 0 \ar[d] \ar[r] & 0 \ar[d] \ar[r] &  G \ar[d] &  \\
          0 \ar[r]  & A' \ar[d] \ar[r] & S \ar[d]^{\sim} \ar[r] & B' \ar[d] \ar[r] & 0 \\
          0 \ar[r]  & A \ar[d] \ar[r] & S \ar[d] \ar[r] & B \ar[d] \ar[r] & 0  \\
          & F \ar[d] \ar[r] & 0 \ar[r] & 0  \\
           & 0  &  &    & 
          } 
\end{equation*}
By the Snake Lemma $G \cong F$. Moreover, since $\langle E_i, E_j \rangle$ abelian, $B' \in \langle E_i, E_j \rangle$.
But now $F \in \langle E_{j_1},\ldots,E_{j_{k-1}} \rangle$. So, $j_{k-1} = i$.
Using again the previous argument, $A, B \in \langle E_i, E_j \rangle$. 

For the second part, $\mc{R}_{E_j} E_i [-1] \in \cat{Q}$ is defined by the exact sequence
$$0 \to \Hom^1 (E_i, E_j) \otimes E_j \to \mc{R}_{E_j} E_i [-1] \to E_i \to 0.$$
If
\begin{equation}\label{12}
0 \to C \to \mc{R}_{E_j} E_i [-1] \to D \to 0
\end{equation}
is a destabilizing sequence and $D \in \langle E_i, E_j \rangle$, then we have a morphism
$$\mc{R}_{E_j} E_i [-1] \to E_{l},$$
with $l = j$ or $l = i$. But this implies $E_{l} \cong E_i$ and so $D \cong E_i$, that is (\ref{12}) is not a destabilizing sequence.        
\end{proof}

It is important to note here that if $\cat{T}$ is the bounded derived category of an abelian category, then, in the assumptions of the previous theorem, $\Tr (E_i , E_j)$ is equivalent to $D (P_k)$, where $k = \dim \Hom^1 (E_i, E_j)$. Indeed, as in the proof of Theorem \ref{thm:Bondal}, we can define a functor $\Yy :  \Tr (E_i , E_j) \to D(P_k)$ as the composition of the inclusion $\Tr (E_i , E_j) \hookrightarrow \cat{T}$ with ${\mathbf R} \Hom (E_i \oplus E_j ,  - )$. The proof that this functor is an equivalence goes in the same way as \cite[Theorem 6.2]{Bond}.

\medskip

We conclude this section by constructing some explicit examples of stability conditions.
Let $\cat{T}$ be, as before, a triangulated category of finite type over $\CC$ and let $\mc{E} = \grf{E_0, \ldots, E_n}$ be a complete exceptional collection on $\cat{T}$. Then the Grothendieck group $K(\cat{T})$ is a free abelian group of finite rank isomorphic to $\ZZ^{\oplus (n+1)}$ generated by $E_0,\ldots, E_n$.
Fix $z_0,\ldots , z_n \in H$. Consider the abelian category $\cat{Q}_p := \langle E_0 [p_0], \ldots, E_n [p_n] \rangle$ of Lemma \ref{lem:tstruct}, for $p_0, \ldots, p_n$ integral numbers such that the exceptional collection $\grf{E_0 [p_0], \ldots, E_n [p_n]}$ is Ext.
Define a stability function $Z_p : K(\cat{Q}_p) \to \CC$ by $Z(E_i [p_i]) = z_i$, for all $i$.
By Remark \ref{rmk:tstruct} this extends to a \emph{unique} stability condition on $\cat{T}$ which is locally-finite. We call the stability conditions constructed in this way \emph{degenerate} if $\rk_{\RR} Z_p = 1$ (seeing $Z_p$ as a map from $K(\cat{T}) \otimes \RR$ to $\CC \cong \RR^2$). Otherwise we call them \emph{non-degenerate}.

Define $\Theta_{\mc{E}}$ as the subset of $\Stab (\cat{T})$ consisting of stability conditions which are, up to the action of $\glpiu$, degenerate or non-degenerate for $\mc{E}$. Then, by Lemma \ref{lem:excobj}, $E_0, \ldots, E_n$ are stable for all stability conditions on $\Theta_{\mc{E}}$; for a degenerate stability condition, they are the only stable objects (up to shifts).

\begin{remark}\label{rmk:stability}
Note that in general $\Theta_{\mc{E}}$ is not the subset of $\Stab (\cat{T})$ consisting of stability conditions in which $E_0, \ldots, E_n$ are stable. Indeed, consider the strong complete exceptional collection $\grf{\mc{O}_{\PP^2}, \mc{O}_{\PP^2} (1), \mc{O}_{\PP^2} (2)}$ on $D (\PP^2)$. Then, by Proposition \ref{pro:mutstab}, in the stability condition $\ss$, defined on the category $\langle \mc{O}_{\PP^2} [2], \mc{O}_{\PP^2} (1) [1], \mc{O}_{\PP^2} (2) \rangle$ by $Z (\ko_{\PP^2} [2]) = Z (\ko_{\PP^2} (1) [1]) = -1$, $Z (\ko_{\PP^2} (2)) = i$, all the objects of the strong complete exceptional collection $\ke := \grf{\ko_{\PP^2} (1), T_{\PP^2}, \ko_{\PP^2} (2)}$ are all stable, but $\ss \notin \Theta_{\ke}$.
\end{remark}

\begin{lemma}\label{lem:theta}
$\Theta_{\mc{E}} \cont \Stab (\cat{T})$ is an open, connected and simply connected $(n+1)$-dimensional submanifold.
\end{lemma}
\begin{proof}
Before proceeding with the proof we need some notations.
Let $\mc{F}_s = \grf{F_0, \ldots, F_s}$, $s>0$, be an exceptional collection.
First of all, define, for $i < j$,
$$k_{i,j}^{\mc{F}_s} :=
\begin{cases}
+\infty, & \text{if}\, \Hom^k (F_i, F_j) = 0,\ \text{for all}\, k;\\
\min \grf{k \, : \, \Hom^k (F_i, F_j) \neq 0}, & \text{otherwise.}
\end{cases}
$$
Then define inductively $\aa_i^{\kf_s} \in \ZZ \cup \grf{+\infty}$ in this way: set $\aa_s^{\kf_s} =0$ and, for $i < s$,
$$\aa_i^{\kf_s} := \us{j>i}{\min} \grf{k_{i,j}^{\kf} + \aa_j^{\kf_s}} - (s-i-1),$$
where the minimum is taken over $\ZZ \cup \grf{+ \infty}$.

Consider $\RR^{n+1}$ with coordinates $\ff_0, \ldots, \ff_n$.
Let $\kf_s := \grf{E_{l_0}, \ldots, E_{l_s}} \cont \grf{E_0, \ldots, E_n}$, $s>0$.
Define $R^{\kf_s}$ as the relation $\ff_{l_0} < \ff_{l_s} + \aa_0^{\kf_s}$.
Finally define
$$\kc_{\ke} := \grf{(m_0, \ldots, m_n, \ff_0, \ldots, \ff_n) \in \RR^{2(n+1)}\, :\begin{array}{l}\bullet\ m_i > 0,\, \text{for all $i$}\\
               															       \bullet\ R^{\kf_s},\, \text{for all } \kf_s \cont \ke,\, s>0 \end{array}}.$$
Clearly $\kc_{\ke}$ is connected and simply connected.

Define a map $\rho : \Theta_{\ke} \to \kc_{\ke}$ by $m_i (\rho(\ss)) := Z(E_i) / |Z(E_i)|$, $\ff_i (\rho(\ss)) := \ff_{\ss} (E_i)$, for $i=0, \ldots, n$, where $\ss=(Z,\kp)$.
We want to prove that this map is an homeomorphism.
By definition of $\Theta_{\ke}$, $\rho$ is injective. Moreover it is straightforward to check that it is also surjective.

Consider the abelian category $\cat{Q}_p$, for $p_0, \ldots, p_n$ integral numbers such that the complete exceptional collection $\grf{E_0 [p_0], \ldots, E_n [p_n]}$ is Ext. Let $\ss_p \in \Theta_{\ke}$ be the stability condition defined by setting $Z_p (E_j [p_j]) = i$, for all $j$ and let $\Gg_p \cont \Stab (\cat{T})$ be the connected component containing $\ss_p$.
First of all notice that the linear subspace $V(\Gg_p)$ of Theorem \ref{thm:bridmain} is all $(K(\cat{T}) \otimes \CC)^\vee$.
Indeed, for all $W \in (K(\cat{T}) \otimes \CC)^\vee$,
\begin{equation*}
\begin{split}
\| W \|_{\ss_p} &= \sup \grf{\frac{| W (E) |}{| Z_p (E) |} \, : \, E \ \mbox{is} \ \ss_p \text{-semistable}} \\
                         &= \sup \grf{\frac{| W (E) |}{| Z_p (E) |} \, : \, E \in \cat{Q}_p} \\
                         &= \sup \grf{\frac{|a_0 W(E_0 [p_0]) + \ldots + a_n W(E_n [p_n])|}{a_0 + \ldots + a_n} \, : \, a_0, \ldots, a_n \geq 0} \\
                         &\leq \sup \grf{\frac{a_0 |W(E_0 [p_0])| + \ldots + a_n |W(E_n [p_n])|}{a_0 + \ldots + a_n} \, : \, a_0, \ldots, a_n \geq 0} \\
                         &< \infty.
\end{split}
\end{equation*}
Hence $\Gg_p$ is a manifold of dimension $(n+1)$. Since the map $\mc{Z}$ of Theorem \ref{thm:bridmain} is a local homeomorphism, to prove the lemma is sufficient to show that $\Theta_{\mc{E}}$ is contained in $\Gg_p$. But, by definition of the generalized metric \eqref{eq:metric}, all stability conditions with heart $\cat{Q}_p$ are in $\Gg_p$ and so is in $\Gg_p$ the open subset $U_p$ consisting of stability conditions which have, up to the action of $\glpiu$, $\cat{Q}_p$ as heart.

Now, let $l$ be an index such that $\Hom^0 (E_s [p_j], E_s [p_l +1]) =0$, for all $s < l$ (for example $l=0$ always works). Then if we set $p'_l := p_l +1$ and $p'_j := p_j$, for $j \neq l$, then the exceptional collection $\grf{E_0 [p'_0], \ldots, E_n [p'_n]}$ is still Ext and $U_p' \cap U_p \neq \emptyset$. Indeed, the stability condition with heart $\cat{Q}_p$ defined by $Z(E_l [p_l]) =i$, $Z(E_j [p_j]) = -1$, for $j \neq l$ is in the $\glpiu$-orbit of the stability condition with heart $\cat{Q}_{p'}$ defined by $Z(E_l [p'_l]) = -1$, $Z(E_j [p'_j]) = i$, for $j \neq l$.

Finally, given two collection of integers $p$ and $q$ such that $\mc{E}_p$ and $\mc{E}_q$ are Ext, it is always possible to find a third collection $r$ such that $r$ is obtained from $p$ and $q$ by successively adding $1$ to some integer corresponding to an index $l$ as before. Hence $\Theta_{\ke} \cont \Gg_p$, for some (all) $p$, as wanted.
\end{proof}

Assume that $\mc{E}$ is regular and that all its iterated mutations are again regular.
Define $\Ss_{\mc{E}}$ as the union of the open subsets $\Theta_{\mc{F}}$ over all iterated mutations $\mc{F}$ of $\mc{E}$. When the triangulated category $\cat{T}$ is \emph{constructible}, i.e. all complete exceptional collections can be obtained, up to shifts, by iterated mutations of a single complete exceptional collection, we simply denote the previous open subset by $\Ss (\cat{T})$.

\begin{corollary}
$\Ss_{\mc{E}} \cont \Stab (\cat{T})$ is an open and connected $(n+1)$-dimensional submanifold.
\end{corollary}
\begin{proof}
It sufficient to show that, for a single mutation $\mc{F}$ of $\mc{E}$, $\Theta_{\mc{F}} \cap \Theta_{\mc{E}}$ is nonempty. We can restrict to consider $\mc{F} = \mc{R}_j \mc{E}$, $j \in \grf{0, \ldots, n-1}$, with $(E_j, E_{j+1})$ not an orthogonal pair.
Fix integers $p_0, \ldots, p_n$ such that  $\Hom^k (E_l [p_l], E_t [p_t]) =0$, for all $k \leq 1$ and all pair $(l, t)$, $l \neq t$ besides $(j,j+1)$, where $\Hom^1 (E_j [p_j], E_{j+1} [p_{j+1}]) \neq 0$. In particular $\grf{E_0 [p_0], \ldots, E_n [p_n]}$ is Ext.
Fix $z_0,\ldots , z_n \in H$ such that $z_l = i$ for $l \neq j, j+1$, $z_j = -1$ and $z_{j+1} = 1+i$ and consider the abelian category $\cat{Q}_p$. Let $\ss$ be the stability condition constructed by these data. Then $\ss \in \Theta_{\mc{F}} \cap \Theta_{\mc{E}}$. Indeed, the exceptional collection
$$\grf{E_0 [p_0], \ldots, E_{j+1} [p_{j+1} +1], \mc{R}_{E_{j+1}} E_j [p_j -1], \ldots, E_n [p_n]}$$
is still Ext and consists of $\ss$-stable objects with phases in the interval $[-1/2, 1/2)$, by Proposition \ref{pro:mutstab}. But then, using the $\glpiu$-action, Lemma \ref{lem:excobj} implies that $\ss$ is in $\Theta_{\mc{F}}$.
\end{proof}

In this way, in the case of projective spaces and del Pezzo surfaces, by Example \ref{ex:PN} and Example \ref{ex:KO}, we constructed open connected subsets of the corresponding stability manifolds. A further study of some topological properties of these subsets is contained in \cite{emolo}.


\section{Examples}\label{sec:examples}

In this section using the description of exceptional objects on the category of representations of quivers without loops given in \cite{CB, Ring}, we study the space $\Stab (P_n)$, which, for $n=2$, is the space of stability conditions on the derived category of $\PP^1$.

Let $P_n$ be the quiver defined in the Example \ref{ex:Pn} and $\mc{Q}_n$ the abelian category of its finite dimensional representations. Since the case $n=0$ is trivial, assume $n >0$.
Set $\grf{S_i}_{i \in \ZZ}$ the family of exceptional objects on $D(P_n)$, where $S_0 [1]$ and $S_1$ are the minimal objects in $\mc{Q}_n$ and the others exceptional objects are defined by
\begin{alignat*}{4}
\ & S_i :=  \mc{L}_{S_{i+1}} S_{i+2}, \quad &  & i < 0,\\
\ & S_i :=  \mc{R}_{S_{i-1}} S_{i-2}, \quad & & i \geq 2.
\end{alignat*}
According to \cite{CB, Ring} these are (up to shifts) the only exceptional objects in $D(P_n)$.
Note that, since $\mc{Q}_n$ is an abelian category of dimension $1$, $S_{\leq 0} [1], S_{>0} \in \mc{Q}_n$ unless $n=1$. In fact the case $n=1$ is also somehow degenerate: indeed in that case three mutations are equal to a shift and so there are effectively only three exceptional objects up to shifts. The main results of this subsection (Lemma \ref{lem:stablepairPn} and Theorem \ref{thm:mainPn}) still hold true for $n=1$; the proofs are a little bit different but easier. Hence we leave them to the reader and in the following we assume $n \geq 2$.

\begin{lemma}\label{lem:helix}
If $i < j$, then
\begin{itemize}
\item $\Hom^k (S_i , S_j) \neq 0$ only if $k = 0$;
\item $\Hom^k (S_j , S_i) \neq 0$ only if $k = 1$.
\end{itemize}
In particular the pair $(S_k, S_{k+1})$ is a complete strong exceptional collection.
\end{lemma}

\begin{lemma}\label{lem:stablepairPn}
In every stability condition on $D(P_n)$ there exists a stable exceptional pair $(E, F)$.
\end{lemma}
\begin{proof}
First of all note that, since $\dim \mc{Q}_n =1$, each object of $D(P_n)$ is isomorphic to a finite direct sum of shifts of objects of $\mc{Q}_n$. So, if an object is stable, some shift of it must belong to $\mc{Q}_n$.

Let $L$ be an exceptional object of $D(P_n)$. We can assume $L$ to be $S_0 [1]$ in $\mc{Q}_n$. Suppose that $L$ is not semistable. Then there exists a destabilizing triangle (the last triangle of the HN filtration)
$$X \to L \os{f}{\to} A  \to X[1],$$
with $A$ semistable and
\begin{equation}\label{13}
\Hom^{\leq 0} (X, A) = 0.
\end{equation}
From the previous remark, we can assume
\begin{align*}
A& = B_0 [0] \oplus B_1 [1],    &   & \ B_i \in \mc{Q}_n\\
f& = f_0 + f_1,                 &   & \ f_0 \in \Hom (L, B_0), \ f_1 \in \Ext^1_{Q_n} (L, B_1).
\end{align*}
Moreover, $f_i =0$ if and only if $B_i = 0$, for $i=0,1$.

But, if $f_0 \neq 0$, then $B_0$ has a direct factor isomorphic to $L$, which is of course not possible.
So, $B_0 = 0$ and the destabilizing triangle is obtained by the extension in $\mc{Q}_n$
\begin{equation}\label{14}
0 \to B_1 \to X \to L \to 0
\end{equation}
corresponding to $f_1$, with the conditions $\Hom (X, B_1) = \Ext^1 (X, B_1) = 0$.
Applying to (\ref{14}) the functor $\Hom ( \bullet , B_1)$ we get $\Ext^1 (B_1, B_1) = \Hom (L, B_1) = 0$.
If $B_1$ is indecomposable, then by \cite{CB, Ring} $B_1$ must be an exceptional object (i.e. $\Hom (B_1, B_1) = \CC$). Otherwise, if $B_1$ is not indecomposable, then again by \cite{CB, Ring}, since there are no orthogonal exceptional pairs\footnote{One should note the different uses of the term ``orthogonal'' here and in \cite{Ring}.}, $B_1$ must be of the form $E^{\oplus i}$ for some exceptional object $E \in \mc{Q}_n$.
But then, by Lemma \ref{lem:helix}, $\Ext^1 (B_1, L) = 0$.
Applying to (\ref{14}) the functor $\Hom ( \bullet , L)$ we get $\Ext^1 (X, L) = 0$. Applying the functor $\Hom (X, \bullet)$ we get $\Ext^1 (X, X) = 0$.
Again $X \cong F^{\oplus j}$, for some exceptional object $F$, and so $(E, F)$ is an exceptional pair.
But also $F$ is semistable. Indeed, suppose it is not; then the HN filtration for $L$ continues
$$R' \to X \to A' \to R'[1],$$
with
\begin{align}\label{16}
\Hom^{\leq 0} & (R', A') = 0, \notag \\
\Hom^{\leq 0} & (A', A) = 0.
\end{align}
Now, proceeding as before, $A'$ and $R'$ are direct sums of exceptional objects. But condition \eqref{16} implies that $A'$ must be $X$, by Lemma \ref{lem:helix}. Hence $F$ is semistable.

To conclude, assume that $L$ is an exceptional object which is semistable but not stable. By repeating the first part of the proof (taking a Jordan-H\"older filtration), we have that there exists an exceptional pair $(E, F)$ of semistable objects such that all the stable factors of $E$ are isomorphic to an object $K_E$ and those of $F$ to another object $K_F$. But then $K_E \cong E$ and $K_F \cong F$. Hence $E$ and $F$ are stable.
\end{proof}

\begin{corollary}
In every stability condition on $D(\PP^1)$ there exists an integer $k$ such that the line bundles $\mc{O} (k)$ and $\mc{O} (k+1)$ are stable.
\end{corollary}

Let $\Theta_k$, $k \in \ZZ$, be the open connected and simply connected subset of $\Stab (P_n) := \Stab (D(P_n))$ defined in the previous section, consisting of stability conditions which are, up to the action of $\glpiu$, degenerate or non-degenerate for the exceptional pair $(S_k, S_{k+1})$.
By Lemma \ref{lem:excobj}, $\Theta_k$ coincides with the subset of $\Stab (P_n)$ consisting of stability conditions in which $S_k$ and $S_{k+1}$ are stable. Moreover, by Lemma \ref{lem:stablepairPn}, $\Stab (P_n)$ is the union over $\ZZ$ of its subsets $\Theta_k$, i.e.\ it coincides with its open subset $\Ss (P_n)$ defined in the previous section. Hence it is connected.
To have a precise description of the topology of $\Stab (P_n)$ we only have to understand how the $\Qq_k$ overlap. The answer is given by the following proposition.

Consider the abelian category $\mc{Q}_n$ and set
\begin{alignat*}{4}
\ & Z^{-1} (S_0 [1]) = -1, \qquad & (& \Rightarrow \ff (S_0 [1]) = 1),\\
\ & Z^{-1} (S_{1}) = 1+i, \qquad & (& \Rightarrow \ff (S_1) = 1/4).
\end{alignat*}
By Remark \ref{rmk:tstruct}, this extends to a unique stability condition $\sigma^{-1} = (Z^{-1}, \mc{P}^{-1})$ on $D(P_n)$.
Consider its $\glpiu$-orbit $O_{-1}$, which is an open subset of $\Stab (P_n)$ homeomorphic to $\glpiu$. Notice that in the case $n=2$, i.e.\ for $\PP^1$, the stability condition induced by $\mu$-stability \cite[Example 5.4]{Brid1} is in $O_{-1}$.

\begin{proposition}\label{pro:stabilityPn}
For all integers $k \neq h$ we have
$$\Theta_k \cap \Theta_h = O_{-1}.$$
\end{proposition}
\begin{proof}
First of all, the fact that $O_{-1} \cont \Theta_k \cap \Theta_h$ is a simple consequence of Proposition \ref{pro:mutstab}.
Let $\ss \in \Theta_k \cap \Theta_h$. Set $\ff_0$ and $\ff_1$ the phases of $S_k [1]$ and $S_{k+1}$ respectively. Then, since $\Hom (S_k, S_{k+1}) \neq 0$, there exists an integer $p \geq -1$ such that
$$\ff_1 -1 < \ff_0 + p \leq \ff_1.$$
But if $p \geq 0$ then there are no stable object in $\ss$ besides $S_k$ and $S_{k+1}$ and so $\ss \notin \Theta_h$. Hence $p=-1$.
If $\ff_0 = \ff_1 +1$, then $S_k$ and $S_{k+1}$ are stable with the same phase, a contradiction. Hence $\ff_1 < \ff_0 < \ff_1 +1$.
We can then act by an element of $\glpiu$ and assuming $S_k [1], S_{k+1} \in \mc{P} ((0, 1])$. By Lemma \ref{lem:excobj} $\langle S_k [1], S_{k+1} \rangle = \mc{P} ((0, 1])$.
So, either $S_0, S_1 \in \mc{P} ((0, 1])$ or $S_0 [1], S_1 [1] \in \mc{P} ((0, 1])$ and, by Proposition \ref{pro:mutstab}, they are both stable. But then acting again with an element of $\glpiu$ we can assume $S_0 [1], S_{1} \in \mc{P} ((0, 1])$ and that the stability function coincides with $Z_{-1}$. But then again by Lemma \ref{lem:excobj} and by Remark \ref{rmk:tstruct}, the resulting stability condition is $\ss^{-1}$.
\end{proof}

\begin{theorem}\label{thm:mainPn}
$\Stab (P_n)$ is a connected and simply connected $2$-dimensional complex manifold.
\end{theorem}
\begin{proof}
By Proposition \ref{pro:stabilityPn}, the simply connected open subsets $\Qq_k$ glue on $O_{-1}$, which is contractible.
The theorem follows from the Seifert-Van Kampen Theorem.
\end{proof}

\begin{corollary}
$\Stab (\PP^1) := \Stab (D(\PP^1))$ is connected and simply connected.
\end{corollary}

The previous corollary has also been obtained independently by Okada \cite{O}. Actually in that paper it is proved a stronger statement: $\Stab (\PP^1) \cong \CC^2$.

Note that the group $\Aut (D(P_n))$ of autoequivalences of $D(P_n)$ acts transitively on the set $\grf{\Qq_k}_{k \in \ZZ}$. Moreover, the subgroup of it which fixes $\Qq_k$ acts trivially on it, up to shifts. But the action of $\Aut (D(P_n))$ on $O_{-1}$ is nontrivial: it is an easy computation to see that
$$O_{-1} / \Aut(D(P_n)) \cong  GL^+ (2, \RR) / G,$$
where $G$ is the subgroup generated by $\begin{pmatrix} 0 & 1\\ -1 & n \end{pmatrix}$.


\section*{Acknowledgements}

It is a pleasure for me to thank Ugo Bruzzo for having introduced me to the problem and for his constant friendly help and encouragement. I also thank Sukhendu Mehrotra, Paolo Stellari and Hokuto Uehara for many comments and suggestions and the Math/Physics research group of the University of Pennsylvania for its warm hospitality during the writing of a first version of this paper. Finally I express my deep gratitude to the referee for the careful reading of the paper and for many comments, suggestions and explanations.


\end{document}